\documentstyle{amsart}

\numberwithin{equation}{section}

\newtheorem{theorem}{Theorem}[section]
\newtheorem{theorema}{Theorem}
\newtheorem{lemma}[theorem]{Lemma}
\newtheorem{corollary}[theorem]{Corollary}
\newtheorem{proposition}[theorem]{Proposition}

\theoremstyle{remark}

\newcommand{\ad}{\operatorname{ad}}
\newcommand{\Ad}{\operatorname{Ad}}

\newcommand{\Fa}{{\frak a}}

\newcommand{\Fg}{{\frak g}}

\newcommand{\Fk}{{\frak k}}
\newcommand{\Fl}{{\frak l}}

\newcommand{\Fp}{{\frak p}}

\newcommand{\g}{{\frak g}}

\newcommand{\cx}{{\Bbb C}}

\begin{document}

\newcommand{\ra}{{\rho^\ast}}
\newcommand{\rb}{{\tilde{\rho}^\ast}}

\title{Prescribing Ricci curvature on complexified symmetric spaces}
\author{Roger Bielawski}

\address{Department of Mathematics, University of Glasgow, Glasgow G12 8QW, UK }

\email{R.Bielawski@@maths.gla.ac.uk}

\maketitle
\thispagestyle{empty}

The aim of this note is to prove
\begin{theorema} Let $G/K$ be an irreducible symmetric space of compact type with $K$ connected and $G$ acting effectively on $G/K$, and let $\rho$ be a real exact $G$-invariant $(1,1)$-form on the complexification $G^\cx/K^\cx$\footnote[1]{The complexification of a compact connected Lie group $G$ is the connected group $G^\cx$ whose  Lie algebra is the complexification of the Lie algebra of $G$ and which satisfies $\pi_1(G^\cx)=\pi_1(G)$.}.
Then there exists a $G$-invariant K\"ahler metric on $G^\cx/K^\cx$ whose Ricci form is $\rho$.
\end{theorema}
{\em Remarks:} 1. The assumptions on $G$ and $K$ are fullfilled for simply connected symmetric spaces. In general, these assumptions guarantee that $G$ is connected (\cite{Hel}, Thm.\ V.4.1), that $G/K$ is globally symmetric (\cite{Hel}, Ex.\ VII.10), and that $K=\Ad_G(K)$ (as $\Ad_G(K)=K/K\cap Z(G)$).\newline
2. The K\"ahler form obtained in Theorem 1 is exact.

\medskip 

The above result has been proved in \cite{Ste} for symmetric spaces of rank $1$ and in \cite{RB} for compact groups, i.e. for the case when $G=K\times K$ and $K$ acts diagonally. For hermitian symmetric spaces and $\rho=0$, Theorem 1 has also been known \cite{BG}.
\par
The proof given here is quite different from that given for group manifolds in \cite{RB}. We show that the complex Monge-Amp\`{e}re equation on $G^\cx/K^\cx$ reduces, for $G$-invariant functions, to a real Monge-Amp\`{e}re equation on the dual symmetric space $G^\ast/K$. We also show that the Monge-Amp\`{e}re operator on non-compact symmetric spaces has a radial part, i.e. it is equal, for $K$-invariant functions, to another Monge-Amp\`{e}re operator on the maximal abelian subspace of $G^\ast/K$. These facts, together with the theorem on $K$-invariant real Monge-Amp\`{e}re equations proved in \cite{RB1}, yield Theorem 1.

\section{Riemannian symmetric spaces of non-compact type}

Here we recall some facts about the geometry of Riemannian symmetric spaces. The standard reference for this section is \cite{Hel}. 
\par
Let $G/K$ be a symmetric space of compact type with $G$ acting effectively and $K$ connected, and let $G^\ast/K$ be its dual\footnote[2]{ Since $K$ is connected and  $K=\Ad_G(K)$ (Remark 1), there exists a $G^\ast$ such that  $(G^\ast,K)$ is a symmetric pair corresponding to $(\Fg^\ast,\Fk)$.}. If $\g$, $\g^\ast$ and $\Fk$ denote the Lie algebras of $G$, $G^\ast$ and $K$, then $\g=\Fk\oplus \Fp$, $\g^\ast=\Fk\oplus i\Fp$, where $[\Fk,\Fp]\subset \Fp$ and  $[\Fp,\Fp]\subset \Fk$. The restriction of the Killing form to $i\Fp$ is positive definite and induces the Riemannian metric of $G^\ast/K$. Moreover, the Riemannian exponential mapping provides a diffeomorphism between $\Fp$ and $G^\ast/K$. This can be viewed as the map:
\begin{equation} p\mapsto e^{ip}K,\label{exp}\end{equation}
where $p\in \Fp$ and $e$ is the group-theoretic exponential map for $G^\ast$. 
Thus we have two $K$-invariant Riemannian metrics on $\Fp\simeq {\Bbb R}^n$: the Euclidean one given by the Killing form, and the negatively curved one given by the diffeomorphism \eqref{exp}.
\par
Let  $\Fa$ be a maximal abelian subspace of $\Fp$ and $\Fl$ its centraliser in $\Fk$. Let $\Sigma$  the set of restricted roots and $\Sigma^+$ the set of restricted positive roots. For each $\alpha\in \Sigma$, let ${\frak p}_\alpha$ (resp. $\Fk_\alpha$) denote the subspace of $\Fp$ (resp. of $\Fk$) where each $(\ad H)^2$, $H\in \Fa$, acts with eigenvalue $\alpha(H)^2$.  We have the direct decompositions
\begin{equation} \Fp=\Fa+\sum_{\alpha\in \Sigma^+}\Fp_{\alpha},\qquad \Fk=\Fl+\sum_{\alpha\in \Sigma^+}\Fk_{\alpha}.\label{decomposition1}\end{equation}
Let $\Fa^+$ be an open Weyl chamber and let $\Fp^\prime$ be the union of $K$-orbits of points in $\Fa^+$. Any $K$ orbit in $\Fp^\prime$ is isomorphic to $K/L$ where the Lie algebra of $L$ is $\Fl$.
Moreover, we have the diffeomorphism:
\begin{equation} \Fa^+\times K/L\rightarrow \Fp^\prime, \qquad (h,k)\mapsto \Ad(k)h.\label{diffeo}\end{equation}
 We now wish to write the two $K$-invariant metrics on $\Fp$ in coordinates given by this diffeomorphism. Let  $\sum dr_i^2$ be the Killing metric on  $\Fa^+$ (the $r_i$ can be viewed as $K$-invariant functions on $\Fp^\prime$). For each  $\Fk_\alpha$, choose a basis  $X_{\alpha,m}$ ($m$ runs from $1$ to twice the multiplicity of $\alpha$) of vectors orthonormal for the Killing form and denote by $\theta_{\alpha,m}$ the corresponding basis of invariant $1$-forms on $K/L$. We have
\begin{proposition} Let $g_0$ be the Euclidean metric on $\Fp$, given by the restriction of the Killing form, and let $g$ be the negatively curved symmetric metric on $\Fp$ given by the  diffeomorphism \eqref{exp}. Then, under the diffeomorphism \eqref{diffeo} the metrics $g_0$ and $g$ can be written in the form
\begin{equation} \sum_i dr_i^2 +\sum_{(\alpha,m)}F(\alpha(r)) \theta_{(\alpha,m)}^2,\label{inv}\end{equation}
where $F(z)=z^2$ for $g_0$, and $F(z)=\sinh^2(z)$ for $g$. \label{metric}\end{proposition}
\begin{pf} Since all these metrics are $K$-invariant, it is enough to compute them at points of $\Fa^+$. Let $H$ be such a point and let $(h,\rho)$, $h\in \Fa$, $\rho\in T_{[1]}K/L$, be a tangent vector to $ \Fa^+\times K/L$ at $(H,[1])$. The vector $\rho$ can be identified with an element of $\sum\Fk_{\alpha}\subset \Fk$.
The corresponding (under \eqref{diffeo}) tangent vector at $H\in \Fp^\prime$ is $h+[\rho,H]$. Computing the Killing form of this vector yields the formula \eqref{inv} with $F(z)=z^2$ for $g_0$. The formula for $g$ follows from a similar computation, using the expression for the differential of the map \eqref{exp} given in \cite{Hel}, Theorem IV.4.1.
\end{pf}

\section{Monge-Amp\`{e}re equation on symmetric spaces}

Let $(M,g)$ be a Riemannian manifold and $u:M\rightarrow {\Bbb R}$ a smooth function. Then the Hessian of $u$ is the symmetric $(0,2)$-tensor $Ddu$ where $D$ is the Levi-Civita connection of $g$. In local coordinates $x_i$, $Ddu$ is represented by the matrix 
\begin{equation} H_{ij}=\frac{\partial^2u}{\partial x_i\partial x_j}-\sum_k \Gamma^k_{ij} \frac{\partial u}{\partial x_k}.\label{Hessian}\end{equation}
We say that the function $u$ is $g$-convex (resp. strictly $g$-convex), if $Ddu$ is non-negative (resp. positive) definite. The Monge-Amp\`{e}re equation on the manifold $(M,g)$ is then
\begin{equation} {\bf M}_g(u):=(\det{g})^{-1}\det Ddu=f\label{MA}\end{equation}
where $f$ is a given function.
\par
Let $(G^\ast/K,g)$ be a symmetric space of non-compact type given by a Cartan decomposition ${\frak g}^\ast={\frak k}+i{\frak p}$. As in the previous section, we identify $M=G^\ast/K$ with ${\frak p}$ and denote by $g_0$ the (flat) metric given by restricting the Killing form to ${\frak p}$. We have:
\begin{theorem} Let $M\simeq {\frak p}$ be a symmetric space of noncompact type and let $u$ be a $K$-invariant (smooth) function on $M$. Then 
\begin{itemize}
\item[(1)] $u$ is $g$-convex if and only if $u$ is  $g_0$-convex (i.e.  convex in the usual sense on ${\frak p}$).
\item[(2)] The following equality of Monge-Amp\`{e}re operators holds: 
$$ {\bf M}_g(u)=F\cdot  {\bf M}_{g_0}(u),$$
where $F:M\rightarrow {\Bbb R}$ is a positive $K$-invariant smooth function depending only on $M$.\end{itemize} \label{solution} \end{theorem}
We have proved in \cite{RB1} a theorem on the existence and regularity of $K$-invariant solutions to Monge-Amp\`{e}re equations on ${\Bbb R}^n$. From this we immediately obtain
\begin{corollary} Let $(G^\ast/K,g)$ be an irreducible symmetric space of noncompact type and let $f$ be a positive smooth $K$-invariant function on $G^\ast/K$. Then the Monge-Amp\`{e}re equation \eqref{MA} has a global smooth $K$-invariant strictly $g$-convex solution.\hfill $\Box$\label{exist}\end{corollary}
We shall now prove Theorem \ref{solution}. In fact we shall prove it in the following, more general situation. Suppose that we are given a $K$-invariant metric on ${\frak p}$ whose pullback under \eqref{diffeo} can be written as (cf. \eqref{inv}): 
\begin{equation} \sum_i dr_i^2 +\sum_{(\alpha,m)}F_{(\alpha,m)}(\alpha(r)) \theta_{(\alpha,m)}^2,\label{inv2}\end{equation}
where $F_{(\alpha,m)}:{\Bbb R}\rightarrow {\Bbb R}$ are smooth functions vanishing at the origin such that $z^{-1}\frac{dF{(\alpha,m)}}{dz}$ is smooth and positive everywhere.  Proposition \ref{metric} implies that the symmetric metric on $G^\ast/K$ is of this form. We claim that Theorem \ref{solution} holds for any metric $g$ of the form \eqref{inv2}. 
\par
In order to simplify the notation, let us write $j$ for the index $(\alpha,m)$ and $\alpha_j$ for $\alpha$ if $j=(\alpha,m)$. The metric $g$ can be now written as
$$ \sum_i dr_i^2 +\sum_{j}F_j(\alpha_j(r)) \theta_{j}^2.$$
We recall the following formula:
$$ 2Ddu=L_{\nabla u}g,$$
where $L$ is the Lie derivative and $\nabla u$ is the gradient of $u$ with respect to the metric $g$.  On the other hand, for any $(0,2)$-tensor $g$ and vector fields $X,Y,Z$, we have:
$$ (L_Xg)(Y,Z)=X. g(Y,Z)- g([X,Y],Z)-g(Y,[X,Z]).$$  
We now compute $L_{\nabla u}g$ on $\Fp^\prime$ with respect to the basis vector fields $\partial/\partial r_i$, $X_j$, where $X_j$ are dual to $\theta_j$. Here $u$ is a $K$-invariant function. The gradient of $u$ is just $\sum\frac{\partial u}{\partial r_i}\frac{\partial }{\partial r_i}$, in particular it is independent of the functions $F_j$.
It follows immediately that $(L_{\nabla u}g)(\partial/\partial r_i,X_j)=0$ and that the matrix  $(L_{\nabla u}g)(X_j,X_k)$ is equal to $\nabla u.g(X_j,X_k)$ and hence it is diagonal with the $(jj)$-entry equal to
$$\nabla u\left(F_j(\alpha_j(r))\right)=
\frac{dF_j}{dz}_{|_{z=\alpha_j(r)}} 
\alpha_j(\nabla_0\bar{u}).$$
Here $\nabla_0\bar{u}=\sum\frac{\partial u}{\partial r_i}\frac{\partial }{\partial r_i}$ is the gradient of $u$ restricted to the Euclidean space $\Fa={\Bbb R}^n$ in coordinates $r_i$, and viewed as a map from ${\Bbb R}^n$ to itself.
\par
Theorem \ref{solution} with the more general metric \eqref{inv2} follows easily with the function $F$ given explicitly by
$$ F=\frac{\prod \alpha_j(r)}{\prod F_j(\alpha_j(r))}\prod \left(\frac{1}{2}\frac{dF_j}{dz}\right)_{z=\alpha_j(r)}.$$
Observe that the assumptions on the $F_j$ guarantee that $F$ extends to a smooth positive function on $\Fp$.

\section{Proof of the Main Theorem}

Let $G/K$ be a locally symmetric space of compact type with $K$ connected.
There is a canonical isomorphism  between $G^\cx/K^\cx$ and $G\times_K {\frak p}$ (i.e. the tangent bundle of $G/K$) given by the map: 
\begin{equation} G\times {\frak p}\rightarrow G^\cx\rightarrow G^\cx/K^\cx, \enskip (g,p)\mapsto ge^{ip}.\label{mostow}\end{equation}
This isomorphism can be viewed in many ways: as an example of Mostow fibration \cite{Mos}, as given via K\"ahler reduction of $G^{\Bbb C}\simeq G\times {\frak g}$ by the group $K$ \cite{HH}, or as given by the adapted complex structure  construction \cite{LS} which provides a canonical diffeomorphism between the tangent bundle of $G/K$ and a complexification of $G/K$. In any case it provides a fibration
\begin{equation} \pi:G^\cx/K^\cx\rightarrow G/K.\label{proj}\end{equation}
The fibers of this projection can be identified with ${\frak p}$ via the map \eqref{mostow}. In particular, the fiber over $[1]$ is given by the $K^\cx$-orbits of elements $e^{ip}, p\in {\frak p}$. We shall relate $G$-invariant plurisubharmonic functions on $G^\cx/K^\cx$ to convex functions on this fiber (see \cite{AL} for a different approach to this).
\par 
For a function $w$ on a complex manifold one defines its Levy form $Lw$ to be the Hermitian $(0,2)$ tensor given in local coordinates as
\begin{equation} \frac{\partial^2 w}{\partial z_k\partial\bar{z}_l}dz_k \otimes d\bar{z}_l.\label{Levy}\end{equation}
This form does not depend on the choice of local coordinates. We shall compute this form for a $G$-invariant function $w$ on  $G^\cx/K^\cx$. It is enough to compute it at points $e^{ip}$, $p\in{\frak p}$. First of all, we choose local holomorphic coordinates at such a point:
\begin{lemma}   In a neighbourhood of a point $e^{ip}$, $p\in{\frak p}$, complex coordinates are provided by the map ${\frak p}^\cx\rightarrow G^\cx\rightarrow G^\cx/K^\cx$, $(a+ib)\mapsto e^{a+ib}e^{ip}$. \label{coor}\end{lemma}
\begin{pf} We have to show that $e^{a+ib}$ does not belong to the isotropy group $e^{ip}K^{\Bbb C}e^{-ip}$ or, equivalently, that for $u\in {\frak p}^\cx$, $\left(\ad e^{-ip}\right)u\not\in {\frak k}^\cx$. We have
\begin{equation} \left(\ad e^{-ip}\right)u=e^{\ad (-ip)}u=\cosh\bigl(\ad (-ip)\bigr)u+ \sinh\bigl(\ad (-ip)\bigr)u,\label{conj}\end{equation}
where the first term of the sum lies in ${\frak p}^\cx$ and the second one in ${\frak k}^\cx$. To show that the first term does not vanish recall that $(\ad (-ip) )^2$ has all eigenvalues nonnegative.\end{pf}

We now have:
\begin{lemma} In the complex coordinates $z=a+ib$ given by the previous lemma, the Levy form \eqref{Levy} of a $G$-invariant function $w$ satisfies the equation:
\begin{equation} \left(\frac{\partial^2 w}{\partial z_k\partial\bar{z}_l}\right)_{\scriptsize\begin{array}{l}a=0\\b=0\end{array}}= \frac{1}{4}\frac{\partial^2}{\partial b_k \partial b_l}w\bigl(e^{ib}e^{ip}\bigr)_{b=0}.\label{key}\end{equation}
\end{lemma}
\begin{pf}
The polar decomposition of $G^\cx$ implies that $e^{a+ib}$ can be uniquely written as $ge^{iy}$, where $g\in G$ and $y\in \g$. Any $G$-invariant function on $G^\cx/K^\cx$  in a neighbourgood of $e^{ip}$ is a function of $y$ only. On the other hand, as $e^{2iy}=\bigl(e^xe^{iy}\bigr)^\ast \bigl(e^xe^{iy})=e^{-a+ib}e^{a+ib}$, it follows from the Campbell-Hausdorff formula that $y=b+[b,a]/2+\text{\sl higher order terms}$. Hence the matrix of second derivatives in \eqref{Levy} at $e^{ip}$ (i.e. at $a=0, b=0$) is the same as the matrix of second derivatives of 
\begin{equation}(a,b)\mapsto e^{(ib+\frac{i}{2}[b,a])}e^{ip} \label{one}\end{equation}
at $a=0, b=0$. We shall now show that for a $G$-invariant function $w$ on $G^\cx/K^\cx$, this matrix of second derivatives is equal to the right-hand side of \eqref{key}.
\par
The Campbell-Hausdorff formula implies that up to order $2$ in $a,b$, we have $e^{(ib+\frac{i}{2}[b,a])}=e^{ib}e^{\frac{i}{2}[b,a])}$. Set $c=[b,a]/2$, which is a point in ${\frak k}$. We are going to show that modulo terms of order $2$ in $c$ (hence of order $4$ in $a,b$), $e^{ic}e^{ip}$ is equal to $e^\rho e^{ip}e^{iq}$, where $\rho\in \g$ and $q\in {\frak k}$ are both linearly dependent on $c$. We note that this proves the lemma, as 
$$e^{ib}e^\rho e^{ip}e^{iq}=e^{\rho}e^{ib+O(3)}e^{ip}e^{iq}= e^{\rho}e^{ib+O(3)}e^{ip}$$
in $G^\cx/K^\cx$, where $O(3)$ denotes terms of order $3$ and higher in
 $a,b$. 
\par
We find $q$ from the equation $\cosh \ad(ip)(q)=c$, which can be solved uniquely as $\cosh \ad(ip)$ is symmetric and positive-definite on ${\frak k}\subset {\frak g}$. We then put $\rho= -i\sinh \ad(ip)(q)$. We observe that $\rho \in \g$ and $e^{\rho-ic}=e^{ip}e^{-iq}e^{-ip}$, thanks to 
\eqref{conj}. Moreover, modulo terms quadratic in $c$, $e^{\rho}=e^{ic}e^{\rho-ic}$ and, consequently:
$$ e^\rho e^{ip}e^{iq}=e^{ic}e^{\rho-ic}e^{ip}e^{iq}= 
e^{ic}\bigl(e^{ip}e^{-iq}e^{-ip}\bigr)e^{ip}e^{iq}=e^{ic}e^{ip},$$
again modulo terms quadratic in $c$. This finishes the proof of the lemma.\end{pf}

According to this lemma, we have to compute $\frac{\partial^2}{\partial b_k \partial b_l}w\bigl(e^{ib}e^{ip}\bigr)_{b=0}$. Now, since $e^{ib}e^{ip}\in G^\ast$, $e^{ib}e^{ip}=ke^{iz}$, where $z=z(b)\in{\frak p}$ and $k\in K$. As $w$ is $G$-invariant, $w\bigl(e^{ib}e^{ip}\bigr)=w(e^{iz})$ and therefore
$$\frac{\partial^2}{\partial b_k\partial b_l}w\bigl(e^{ib}e^{ip}\bigr)_{b=0}=
\frac{\partial^2}{\partial b_k \partial b_l}w\bigl(e^{iz(b)}\bigr)_{b=0}.$$
Thus we compute the matrix of second derivatives of a function defined on $\exp(i{\frak p})$ in the coordinates given by $b\mapsto e^{ib}e^{ip}\mapsto e^{iz(b)}$. These, however, are the geodesic coordinates at the point $e^{ip}$ in the symmetric space $K\backslash G^\ast$ (being translations of geodesics at $[1]$), and hence the matrix of second derivatives in these coordinates is equal to the Riemannian Hessian \eqref{Hessian} for the symmetric metric on $K\backslash G^\ast$.
\newline 
We obtain
\begin{theorem} Let $w$ be a smooth $G$-invariant function on $X=G^\cx/K^\cx$ and let $\bar{w}$ be its restriction to the fiber $S=\exp(i{\frak p})$ of \eqref{proj} over $[1]$. Let $g$ denote the symmetric metric on $S\simeq K\backslash G^\ast$. Then $w$ is (strictly) plurisubharmonic if and only if $\bar{w}$ is (strictly) $g$-convex. Moreover, the following equality holds:
$$\partial \bar{\partial}\log\det Lw=  \partial \bar{\partial} \log \widehat{{\bf M}_g(\bar{w})},$$
where $\hat{u}:X\rightarrow {\Bbb R}$ is a $G$-invariant function such that $\bar{\hat{u}}$ is a given $K$-invariant function $u$ on $S$.\hfill $\Box$
\label{last}\end{theorem}

We are now ready to prove Theorem 1. Recall that $X=G^\cx/K^\cx$ is a Stein manifold and so if $\rho$ is an exact $(1,1)$ form on $X$, then $\rho=-i\partial\bar{\partial}h$ for some function $h$. If $\rho$ is $G$-invariant, then we can assume that $h$ is $G$-invariant. We can restrict $h$ to the fiber $S$ defined in the last theorem and thanks to Corollary \ref{exist} we can find a strictly $g$-convex $K$-invariant smooth solution $\bar{u}$ to the equation \eqref{MA} with $f=e^h$, where the metric $g$ is the symmetric metric on $S\simeq K\backslash G^\ast$. We can extend this solution via $G$-action to a $G$-invariant function $u$ on $X$. Theorem \ref{last} implies now that $u$ is strictly plurisubharmonic and that the Ricci form of the K\"ahler metric with potential $u$ is $\rho$. 

\medskip

{\em Acknowledgement}. This work has been supported by an Advanced Research Fellowship from the Engineering and Physical Sciences Research Council of Great Britain.

\end{document}